\documentclass[a4paper,10pt,oneside]{amsart}

\usepackage{amsmath,amssymb, amsthm}
\usepackage[english]{babel}
\usepackage{graphicx}
\usepackage{xcolor}
\usepackage{multirow}
\usepackage{url}
\usepackage{hyperref}
\usepackage[all]{xy}
\hypersetup{
    colorlinks,
    linkcolor={red!50!black},
    citecolor={blue!50!black},
    urlcolor={blue!80!black},
    breaklinks=true
}

\newtheorem{lemma}{Lemma}[section]
\newtheorem{proposition}[lemma]{Proposition}
\newtheorem{theorem}[lemma]{Theorem}
\newtheorem{corollary}[lemma]{Corollary}

\newtheorem{remark}[lemma]{Remark}

\newtheorem{definition}[lemma]{Definition}

\newcommand{\p}{\mathbb{P}}

\newcommand{\R}{\mathbb{R}}
\newcommand{\N}{\mathbb{N}}
\newcommand{\Z}{\mathbb{Z}}
\newcommand{\F}{\mathbb{F}}

\renewcommand{\tilde}{\widetilde}
\newcommand{\mcal}[1]{\mathcal{#1}}

\newcommand{\Gal}{\operatorname{Gal}}
\newcommand{\GS}{\operatorname{GS}}

\newcommand{\Disc}{\operatorname{Disc}}
\newcommand{\Hom}{\operatorname{Hom}}

\begin{document}

\title{Probabilities of incidence between lines and a plane curve over finite fields}
\thanks{This article was published in Finite Fields and Their Applications, Vol 61, January 2020, Authors: Mehdi Makhul, Josef Schicho, and Matteo Gallet, Title: ``Probabilities of incidence between lines and a plane curve over finite fields'', Pages 101582, \href{https://doi.org/10.1016/j.ffa.2019.101582}{DOI:10.1016/j.ffa.2019.101582}, Copyright Elsevier (2019).\\
\copyright\ 2019.  This manuscript version is made available under the CC-BY-NC-ND 4.0 license \url{http://creativecommons.org/licenses/by-nc-nd/4.0/}.}

\author{Mehdi Makhul${}^{\ast}$}
\address[MM]{Johann Radon Institute for Computational and Applied 
	Mathematics (RICAM), Austrian Academy of Sciences, Linz\\
	Research Institute for Symbolic Computation (RISC), 
	Johannes Kepler University, Linz}
\email{mmakhul@risc.jku.at}

\author{Josef Schicho${}^{\diamond}$}
\address[JS]{Research Institute for Symbolic Computation (RISC), 
	Johannes Kepler University, Linz}
\email{jschicho@risc.jku.at}

\author{Matteo Gallet${}^{\bullet}$}
\address[MG]{International School for Advanced Studies/Scuola Internazionale Superiore di Studi Avanzati (ISAS/SISSA),
Via Bonomea 265, 34136 Trieste, Italy}
\email{mgallet@sissa.it}

\thanks{${}^{\ast}$ Supported by the Austrian Science Fund (FWF): W1214-N15 
Project DK9, and P30405-N32.} 

\thanks{${}^{\diamond}$ Supported by the Austrian Science Fund (FWF): W1214-N15 
Project DK9.}

\thanks{${}^{\bullet}$ Supported by the Austrian Science Fund (FWF): W1214-N15 
Project DK9, P26607, P25652, P31061 and Erwin Schr\"odinger Fellowship J4253.}


\begin{abstract}
We study the probability for a random line to intersect a given plane curve, 
over a finite field, in a given number of points over the same 
field. In particular, we focus on the limits of these probabilities under 
successive finite field extensions. Supposing absolute irreducibility for the 
curve, we show how a variant of the Chebotarev density theorem for function fields 
can be used to prove the existence of these limits, and to compute them under a 
mildly stronger condition, known as simple tangency. Partial results have 
already appeared in the literature, and we propose this work as an introduction 
to the use of the Chebotarev theorem in the context of incidence geometry. 
Finally, Veronese maps allow us to compute similar probabilities of 
intersection between a given curve and random curves of given degree. 
\end{abstract}

\maketitle 

\section{Introduction}
\label{introduction}

What is the probability that a random line in the (affine or projective) plane 
intersects a curve of given degree in a given number of points? More precisely, 
what happens if we consider a finite field with $q$ elements as base field, and 
then we ask the same question for a field with $q^2, q^3, \ldots, q^N$ elements, 
analyzing how these probabilities behave as $N$ goes to infinity? In this work 
we investigate this problem by means of algebro-geometric techniques. Recently, 
the interplay between combinatorial problems and algebraic techniques has 
become more and more common, and has been revealing to be extremely fruitful. 
Here we refer in particular to the area called \emph{combinatorial geometry}, 
which is described in the abstract of~\cite{Tao2014} as the area dealing with
``\emph{the possible range of behaviours of arbitrary finite collections of geometric 
objects such as points, lines, or circles with respect to geometric operations such as incidence or distance}''.
As Tao points out in~\cite{Tao2014}, in the last decade algebraic geometry and algebraic topology helped to unriddle
several important questions and conjectures in this area. Amongst the most 
prominent of such problems, one can mention the \emph{distinct 
distance problem} (see~\cite{Guth2015}), the \emph{Kakeya problem over finite 
fields} (see \cite{Dvir2009} and later improvements in~\cite{Saraf2008} 
and~\cite{Dvir2013}) and the \emph{Dirac-Motzkin conjecture} 
(see~\cite{Green2013}). For a nice survey about these topics, 
see~\cite{Tao2014}.

A motivation for the problem we investigate in our paper comes from the famous 
\emph{Sylvester-Gallai theorem}. Sylvester posed it as a question 
in~\cite{Sylvester1893}, which was raised again by Erd\"os in~\cite{Erdoes1943} 
and later solved by Melchior (see~\cite{Melchior1941}) and Gallai 
(see~\cite{Gallai1944}). Given a set of points in the affine plane, a line is 
called \emph{ordinary} if it passes through exactly two of them.
\begin{theorem}[Sylvester-Gallai]
Suppose that $P$ is a finite set of points in the real plane, not all on a 
line. Then $P$ admits an ordinary line. 
\end{theorem}
This theorem is clearly false if we work over finite fields, since in this case 
we can pick~$P$ to be the the whole plane. Moreover, the theorem is false also 
on the complex plane: in fact, in \cite{Naimpally1966} Serre observed that the 
$9$ inflection points of a cubic curve do not satisfy the requirements of 
Silvester-Gallai theorem; also, the so-called \emph{$3$-nets} provide other 
counterexamples, see \cite{Korchmaros2015, Yuzvinsky2004}.

In the same circle of ideas, in \cite{Solymosi2013}, Solymosi considered the 
following situation. Given a set~$P$ of points in the plane, a line is called   
$k$-\emph{rich}, if it contains precisely $k$ points of $P$. For example, a 
$2$-rich line is an ordinary line. Then, Solymosi's theorem reads as:
\begin{theorem}[Solymosi]
For any $k \ge 4$, there is a positive integer $n_0$ such that for $n > 
n_0$ there exists $P \subseteq \R^2$ such that there are at least 
$n^{2-\frac{c}{\sqrt{\log n}}}$ $k$-rich lines, but no $k+1$-rich lines. Here, 
$c=2\log(4k+1)$.
\end{theorem}

A recent outstanding result of Green and Tao (see~\cite{Green2013}) gives an 
almost complete description of the structure of sets with few ordinary 
lines in the real plane. In the same paper, the authors also proved the 
Dirac-Motzkin conjecture and a less known problem, referred in the literature 
as the \emph{orchard problem}. 

Our work is inspired by these results. Here, we consider an 
algebraic plane curve~$C$ of degree~$d$ over a finite field~$\F_q$ with $q$ 
elements, where $q$ is a prime power, namely the set of points in the projective 
plane~$\p^2_{\F_q}$ that are zeros of a homogeneous trivariate polynomial of 
degree~$d$. Given such a curve, we can define the probability for a line 
in~$\p^2_{\F_q}$ to intersect it in exactly $k$ points. Notice that here we 
consider the mere set-theoretic intersection: no multiplicities are taken into 
account. We can then consider the same kind of probability, keeping the same 
curve~$C$ --- namely, the same trivariate polynomial --- but changing the base 
field from~$\F_q$ to~$\F_{q^2}$, $\F_{q^3}$ and so on. In this way, for every $N 
\in \N$ we define the numbers~$p_k^N(C)$, namely the probability for a line 
in~$\p^2_{\F_{q^N}}$ to intersect~$C$ in exactly~$k$ points. If the limit as~$N$ 
goes to infinity of the sequence $\bigl( p_k^N(C) \bigr)_{N \in \N}$ exists, we 
denote this number by~$p_k(C)$. The main tool we use to compute these numbers 
when the curve $C$ is absolutely irreducible and with \emph{simple tangency} is 
an effective version of the \emph{Chebotarev theorem} for function fields. 
Here, by \emph{absolutely irreducible} we mean that the curve is irreducible 
over the algebraic closure of its field of definition. By asking that the curve 
has \emph{simple tangency} we require that there exists a line whose 
intersection with~$C$ consists of simple intersections except for one, which is 
a double intersection. These are the main results of our paper:
\begin{theorem}
\label{thm:existence}
Let $C$ be an absolutely irreducible plane algebraic curve of degree~$d$ 
over~$\F_q$, where $q$ is a prime power. Then the numbers $\{ p_k(C) \}$ are 
well-defined, namely the corresponding limits exist.
\end{theorem}
\begin{theorem}
\label{thm:formula}
Let $C$ be an absolutely irreducible plane algebraic curve of degree~$d$ 
over~$\F_q$, where $q$ is a prime power. Suppose that $C$ has simple tangency. 
Then for every $k \in \{0, \dotsc, d\}$ we have 
\[
 p_k(C) = \sum_{s=k}^d \frac{(-1)^{k+s}}{s!} \binom{s}{k}.
\]
In particular, $p_{d-1}(C) = 0$ and $p_d(C) = 1 / d!$.
\end{theorem}

We approached this problem using Galois theory techniques; during a revision of 
our work, we have been informed\footnote{We thank two anonymous referees for 
pointing us to the relevant literature.} that some of the questions investigated 
in this paper (or similar ones) have already appeared in the literature, though 
expressed in a different language and with different purposes 
(see~\cite{Bary-Soroker2012} and~\cite{Diem2012}). Both the two cited paper use 
the Chebotarev theorem for function field as a key ingredient. After studying the
Chebotarev theorem, we realized that we could use it to provide a much shorter 
proof for our result than the one we initially used, and that our initial 
approach, although we were not aware of that, did not differ too much from the 
techniques that lead to the Chebotarev theorem. However, we still think that our 
initial approach could be of interest for researchers in discrete and 
combinatorial geometry. In fact, although it provides less information than the
Chebotarev theorem, it can serve  as an introduction to this technique because 
of its self-containedness and of the avoidance of technical aspects that are 
present in other works. Because of this, in the initial part of this paper we 
report our initial approach to the problem, and then we explain how to use the
Chebotarev theorem to obtain Theorem~\ref{thm:formula}. After that, we show how 
the same technique provides a formula for the probabilities of intersection 
between a given plane curve of degree~$d$ and a random plane curve of degree~$e$ 
(Proposition~\ref{prop:higher_formulas}). Claus Diem then pointed out to us that the material we present 
has essentially already appeared in SGA1~\cite{SGA2003} by Grothendieck, but it is ``a bit 
hidden'', as he said; he suggested another way of presenting 
the material, which we found better than the one we used, and we adopted this choice of exposition\footnote{We 
thank Claus Diem for his careful description of the alternative way of 
presenting the material.}.

We briefly summarize how the problem we investigate is discussed in the 
aforementioned literature. In~\cite{Bary-Soroker2012}, the focus is a variant of 
the so-called \emph{Bateman-Horn conjecture} for polynomial rings of finite 
fields. The original Bateman-Horn conjecture concerns the frequency of prime 
numbers among the values of a system of polynomials at integer numbers. One of 
its consequences is \emph{Schinzel conjecture}, which asks whether, given 
polynomials $f_1, \dotsc, f_r \in \Z[x]$, then for infinitely many $n \in \Z$ we 
have that $f_1(n), \dotsc, f_r(n)$ are all prime. Bary-Soroker and Jarden 
consider the situation in which $\Z$ is replaced by~$\F_q[t]$ for some prime 
power~$q$. More precisely, given polynomials $f_1, \dotsc, f_r \in \F_q[t][x]$, 
they want to compute the number of polynomials $g \in \F_q[t]$ such that $f_1 
\bigl( t,g(t) \bigr), \dotsc, f_r \bigl( t,g(t) \bigr)$ are irreducible. In 
particular, they focus on the case when $g$ is linear, namely on the computation 
of the pairs $(a_1, a_2) \in \F_q^2$ such that $f_1(t, a_1t + a_2), \dotsc,  
f_r(t, a_1t + a_2)$ are irreducible. In our language, this is the number of 
lines in the plane such that the polynomial obtained by restricting a plane 
curve on such a line is irreducible. The authors improve a result by Bender and 
Wittenberg (see \cite[Theorem~1.1 and Proposition~4.1]{Bender2005}) and show 
that this number goes as~$q^2 / d$. To prove this, they make use of an effective 
version of the Chebotarev density theorem (see the appendix of 
\cite{Andrade2015}). The number computed by Bary-Soroker and Jarden is similar 
to the quantity~$p_0$ that we define, though it is not the same, since it can 
happen that a line does not intersect a curve at any point over~$\F_q$, but the 
polynomial given by the restriction of the curve to the line can be reducible. 
Also the behaviour as $d \to \infty$ of these two quantities is different: the 
one by Bary-Soroker and Jarden goes to zero, while $p_0$ tends to~$1/e$.

In \cite{Diem2012}, the author focuses on the complexity of computation of the 
so-called \emph{discrete logarithm} in the group of divisors of degree~$0$ of a 
nonsingular curve. Given two elements $a$ and $b$ in a group~$G$, the discrete 
logarithm~$\log_b a$ is an integer~$k$ such that $b^k = a$. On a smooth 
curve~$C$, one can consider formal integer sums of points of~$C$, and define an 
equivalence relation on them in order to obtain the class group of~$C$. One can 
therefore try to compute discrete logarithms in the class group of a curve, and 
in particular for those formal integer sums of points whose coefficients add up 
to zero, namely the ones of degree~$0$; this has important applications in 
cryptography. In \cite[Theorem~2]{Diem2012}, Diem proves that computing the 
discrete logarithm has an expected time of $\tilde{O}(q^{2- \frac{2}{d-2}})$ 
for those curves over~$\F_q$ that admit a birational plane model~$D$ of 
degree~$d$ such that there exists a line in the plane intersecting~$D$ in $d$ 
distinct points over~$\F_q$. Then the author computes the number of lines in 
the plane intersecting~$D$ in exactly~$d$ points over $\F_q$ (see 
\cite[Theorem~3]{Diem2012}), namely the quantity~$p_d(D)$ in our language. As 
in the previous paper, this is done using an effective version of the Chebotarev 
density theorem (see \cite{KumarMurty1994}).

Recently, a new paper~\cite{Entin2018} appeared dealing with the same problem 
we investigate in our work, but allowing the given curve to be constituted of 
several irreducible components. We have been informed by Kaloyan Slavov that 
also Birch and Swinnerton-Dyer investigated in~\cite{Birch1959} this topic, 
providing a formula for the quantity $p_1(C) + \dotsb + p_d(C)$ in our language. 
We thank him for pointing out to us this reference, and for useful suggestions.

The rest of the paper is structured as follows. Section~\ref{preliminaries} 
introduces some preliminary results, namely the Lang-Weil bound for the number 
of points of a variety over a finite field 
(Subsection~\ref{preliminaries:lang_weil}) and some known facts about Galois 
groups of plane curves (Subsection~\ref{preliminaries:galois}). 
Sections~\ref{galois} and~\ref{probabilities} present our initial approach to 
the problem, which provides less information than the one obtained via the
Chebotarev theorem, but uses more elementary tools (namely, some basic facts 
about \'etale maps). Section~\ref{chebotarev} shows how to use the Chebotarev 
density theorem in order to prove Theorems~\ref{thm:existence} 
and~\ref{thm:formula} and Proposition~\ref{prop:higher_formulas}.

\section{Preliminaries}
\label{preliminaries}

\subsection{Lang-Weil bound}
\label{preliminaries:lang_weil}

One of the main tools we use in our work is the so-called \emph{Lang-Weil 
bound} for the number of points of a variety over a finite field (see 
\cite[Theorem~1]{Lang1954}). For a nice exposition of this result, see Terence 
Tao's blog\footnote{\href{
		https://terrytao.wordpress.com/2012/08/31/the-lang-weil-bound/}{
		https://terrytao.wordpress.com/2012/08/31/the-lang-weil-bound/}}. 
Let $\F$ be a field and consider an affine algebraic variety~$V$ 
over~$\F$. This means that we are given finitely many polynomials $P_1, \dotsc, P_r \in \F[x_1, \dotsc, x_n]$, which generate the so-called ideal of~$V$, denoted $I(V)$.
For any extension of fields $\F \subseteq \mathbb{K}$, we denote 
by~$V(\mathbb{K})$ the set of common zeros in~$\mathbb{K}^n$ of the polynomials 
in the ideal~$I(V)$, considered now as an ideal in $\mathbb{K}[x_1, \dotsc, 
x_n]$. One says that a variety~$V$ is \emph{irreducible} if 
$I(V)$ is prime in~$\F[x_1, \dotsc, x_n]$. For our considerations we will need a 
stronger notion of irreducibility, which we introduce in the following 
definition.

\begin{definition}
 We say that an affine variety~$V$ over a field~$\F$ is \emph{absolutely 
irreducible} if the ideal~$I(V)$ is prime in $\overline{\F}[x_1, \dotsc, x_n]$, 
where $\overline{\F}$ is an algebraic closure of~$\F$.
\end{definition}

\begin{definition}
 We say that an affine variety $V \subseteq \F^n$ defined by 
polynomials $P_1, \dotsc, P_r$ has \emph{complexity}~$M$ if $n, r \leq M$ and 
$\deg(P_i) \leq M$ for all $i \in \{1, \dotsc r\}$.
\end{definition}

\begin{theorem}[Lang-Weil bound]
\label{thm:lang_weil}
Let $V$ be an absolutely irreducible variety over a finite field $\F$ of 
complexity at most~$M$. Then 
\[
 |V(\mathbb{F})| = 
 \bigl( 
  1 + O_M(|\mathbb{F}|^{-\frac{1}{2}}) 
 \bigr) \thinspace 
 |\mathbb{F}|^{\dim(V)} \,.
\]
\end{theorem}
By writing $O_M(|\F|^{-\frac{1}{2}})$ we mean that there exists a 
nonnegative constant~$\delta_M$ depending on~$M$, but not on~$V$, such that 
\[
 \bigl( 1-\delta_M |\F|^{-\frac{1}{2}}) \bigr) \thinspace |\F|^{\dim(V)} 
 \leq 
 |V(\F)| 
 \leq 
 \bigl( 1+\delta_M |\F|^{-\frac{1}{2}}) \bigr) \thinspace |\F|^{\dim(V)} \,.
\]
Using an inclusion-exclusion argument, one obtains by induction on the 
dimension:
\begin{corollary}
\label{cor:lang_weil}
 Let $V$ be a variety over a finite field $\F$ of complexity at most~$M$. Then
\[
 |V(\mathbb{F})| = 
 \bigl( 
  c + O_M(|\mathbb{F}|^{-\frac{1}{2}}) 
 \bigr) \thinspace 
 |\mathbb{F}|^{\dim(V)} \,,
\]
where $c$ is the number of irreducible components of~$V$ that are absolutely irreducible.
\end{corollary}

All the considerations and results we stated so far hold also for 
projective varieties over finite fields. By a \emph{projective variety} over a 
field~$\F$ we mean a variety in the projective space~$\p^n_{\F}$ given by finitely many \emph{homogeneous} polynomials $P_1, \dotsc, P_r \in \F[x_0, \dotsc, x_n]$. From now on, all the varieties we consider are 
projective, or are open subsets of projective varieties.

\subsection{Galois group of a plane curve}
\label{preliminaries:galois}

The aim of this section is to recall a construction (see~\cite{Rathmann1987}) 
which associates a Galois group to a plane algebraic curve. We will see in the 
following sections that this group determines the irreducibility of certain 
surfaces; this will be the key to derive a formula for the probabilities 
we are interested in.

Let $q$ be a prime power, namely $q = p^r$ for some prime number~$p$. We 
denote by~$\F_q$ the finite field with $q$ elements. Let $C$ be an absolutely 
irreducible algebraic curve in~$\p^2_{\F_q}$. 
Define $X_1$ to be the unique subvariety of $C \times \check{\p}_{\F_q}^2$ --- 
here $\check{(\cdot)}$ denotes the \emph{dual projective plane} --- such that,
for every extension field~$K$ of~$\F_q$,
\[
 X_1(K) = 
 \bigl\{ 
  (w,[\ell]) \in C(K) \times \check{\p}^2 (K) :  w \in \ell 
 \bigr\}
 \quad \textrm{and} \quad
 X_0 := \check{\p}^2_{\F_q}.
\]
For a line $\ell \subseteq \p^2_{K}$, we 
write~$[\ell]$ for the corresponding point in~$\check{\p}^2 (K)$. The 
correspondence is given by
\[
 \check{\p}^2 (K) \ni (a:b:c) 
 \quad \longleftrightarrow \quad
 \bigl\{ (x:y:z) \in \p^2 (K) \, : \, ax + by + cz = 0 \bigr\}.
\]

\begin{definition}
\label{def:projection}
Using the notation we have already introduced, we define the map $\pi 
\colon X_1 \longrightarrow X_0$ to be the projection onto the second component.
\end{definition}

Since $X_0$ is irreducible, we can define its \emph{function field}, 
denoted~$K(X_0)$. This is the field of equivalence classes of morphisms $\varphi 
\colon U \longrightarrow \mathbb{A}_{\F_q}^1$, where $U$ is any (Zariski) open subset of~$X_0$; 
two morphisms are considered equivalent if they agree on a non-empty open 
subset. Consider the projection $\rho \colon X_1 \longrightarrow C$ on the first 
component: its fibers are lines in the dual projective space. Hence all these fibers are irreducible 
varieties of the same dimension. This implies that $X_1$ is irreducible by 
\cite[Exercise 11.4.C]{Vakil2017}; its function field is denoted~$K(X_1)$. 

\begin{lemma}[{see \cite[Definition~1.3]{Rathmann1987}}]
\label{lemma:projection_quasi_finite}
	The projection $\pi \colon X_1 \longrightarrow X_0$ is a quasi-finite 
dominant separable morphism of degree~$d$.
\end{lemma}

Because of Lemma~\ref{lemma:projection_quasi_finite}, the induced map $\pi^{*} 
\colon K(X_0) \longrightarrow K(X_1)$ between fields of rational functions 
realizes~$K(X_1)$ as a finite separable extension of~$K(X_0)$ of degree~$d$. By 
the primitive element theorem, the field~$K(X_1)$ is generated over~$K(X_0)$ by 
a single rational function~$h \in K(X_1)$ satisfying $P(h)=0$ for an 
irreducible monic polynomial~$P$ over~$K(X_0)$ of degree~$d$.

\begin{definition}[Galois group, {see \cite[Definition~1.3]{Rathmann1987}}]
\label{def:galois_curve}
Using the notation just introduced, we define the \emph{Galois group}~$\Gal(C)$ 
of~$C$ to be the Galois group of a splitting field of the polynomial~$P$ 
over~$K(X_0)$. In other words, $\Gal(C)$ is the Galois group of a Galois 
closure (see \cite[Remark~4.77]{Rowen2006}) of the field extension $K(X_0) 
\hookrightarrow K(X_1)$. The 
group~$\Gal(C)$ is independent of the choice of~$h$ and it can be regarded as a 
subgroup of the permutation group~$S_d$ of the roots of~$P$.   
\end{definition}

\begin{definition}[Simple tangency]
\label{def:simple_tangency}
Let $C$ be an absolutely irreducible curve of degree~$d$ 
in~$\p^2_{\F_{q}}$. We say that $C$ 
has \emph{simple tangency} if there exists a line~$\ell \subseteq 
\p^2_{\overline{\F}_{q}}$ intersecting~$C$ in $d-1$ smooth points of~$C$ such 
that $\ell$ intersects~$C$ transversely at $d-2$ points and has intersection 
multiplicity~$2$ at the remaining point.
\end{definition}

\begin{remark}
 A general curve $C \subseteq \p^2_{\F_q}$ of degree~$d$ has simple 
tangency. In fact, notice that having simple tangency is an open condition, 
therefore it is enough to exhibit a single example in order to obtain the 
claim. To do that, consider the curve of equation
\[
 x^2 \, P(x,y) + z \, Q(x,y,z) = 0,
\]
where $P$ is a homogeneous polynomial with $d-2$ distinct roots in 
$\overline{\F}_q$ and $Q$ is a homogeneous polynomial of degree~$d-1$.
\end{remark}

\begin{proposition}[{\cite[Proposition~2.1]{Rathmann1987}}]
\label{prop:symmetric_galois}
Let $C \subseteq \p^2_{\F_{q}}$ be an absolutely irreducible plane curve of 
degree~$d$ with simple tangency. Then the Galois group $\Gal(C)$ of~$C$ is the 
whole symmetric group~$S_d$. 
\end{proposition}

Claus Diem pointed out to us that in the original proof of 
Proposition~\ref{prop:symmetric_galois} it is written that ``\emph{For $k = 2$ 
the variety~$U_2$ is a $\p^{n-2}$-bundle over~$C$ and therefore irreducible}''. 
He explained us that this is impossible for dimension reasons. A first 
attempt for a fix would be to replace~$C$ by~$C \times C$. It turns out that 
then the fibers are not (always) proper, and then one cannot conclude that $U_2$ 
is irreducible. A correct argument has already been given by Ballico and Hafez 
in~\cite{Ballico1986}.

\section{Galois theory for \'etale maps}
\label{galois}

In this section we associate a Galois group to a morphism (satisfying 
certain conditions) between two irreducible smooth varieties. We show  that this concept admits a geometric 
counterpart, and we use this characterization in the next section. As we 
pointed out in the Introduction, the results of this and the following section 
are subsumed by the ones of Section~\ref{chebotarev}. Nevertheless, we believe 
that the approach presented in these sections can be useful to help 
understanding the setting that is used also in the Chebotarev theorem to solve this 
kind of problems. Claus Diem pointed out that the material in this 
section is essentially already present in SGA1 \cite{SGA2003}; moreover, he
suggested us a clearer and shorter way to present the material about Galois 
closures of \'etale maps. We follow his suggestions, and we thank him 
for sharing with us this material.

For technical reasons, we develop the theory for a special class of morphisms, 
namely the one of \emph{\'etale} maps. They model, in the algebraic setting, the 
notion of ``local isomorphism'' for the analytic topology. Recall that, in 
differential geometry, a smooth map between two smooth manifolds is a 
\emph{local diffeomorphism} if it induces an isomorphism at the level of tangent 
spaces. For an affine variety~$X$ cut out by polynomials $P_1, \dotsc, P_r$, one 
defines the \emph{tangent cone} $C_{O}(X)$ of~$X$ at the origin as the variety 
defined by the homogeneous parts of minimal degree of each of the 
polynomials~$P_1, \dotsc, P_r$; the tangent cone at any other point is obtained 
by translating it to the origin and by applying the previous definition. The 
tangent cone plays for \'etale morphisms the role played by the tangent space 
for local diffeomorphisms. A morphism $f \colon X \longrightarrow Y$ between varieties over an algebraically closed field is \emph{\'etale at a point} $x \in X$ if it induces an isomorphism between the tangent cones $C_{x}(X)$ and $C_{f(x)}(Y)$. A map is called \emph{\'etale} if it is \'etale at every point (see~\cite[Chapter~2]{Milne2013}). For more general varieties, one adopts the definition of an \'etale map as a morphism which is flat and unramified (see \cite[Chapter~1]{Milne1980}).

We are going to define a notion of Galois closure for \'etale maps.

\begin{remark}
	Consider a separable extension of fields $K \subseteq L$. We can define the 
\emph{Galois closure} of this extension as the minimal extension~$M$ of~$L$  
which is Galois over~$K$. In the language of schemes, the extension $K \subseteq L$ 
corresponds to a connected \'etale map. The two varieties of this map have 
each a single point, but the structure of schemes still allows to encode the 
field extension. We can hence consider the classical notion of Galois closure 
for field extensions as the ``toy'' case of the notion of Galois closure of 
\'etale maps.
\end{remark}

We mimic the classical notion of Galois closure for field extensions in the 
context of maps. The Galois closure of a map, then, is defined as a map 
satisfying a universal property similar to the one satisfied by the Galois 
closure of a field extension.

\begin{definition}
	Let $g \colon Z \longrightarrow X$ be a connected \'etale map, namely 
	both~$Z$ and~$X$ are connected. We say that $Z$ is \emph{Galois} over~$X$ 
	if the group of automorphisms of~$Z$ is transitive on the geometric fibers of~$g$. 
\end{definition}

\begin{definition}
\label{def:galois_closure}
	Let $f \colon X \longrightarrow Y$ be a connected \'etale morphism. A 
\emph{Galois closure} of~$f$ is an \'etale morphism $Z \longrightarrow X$ such 
that $Z$ is Galois over~$Y$ and such that $Z \longrightarrow X$ is minimal under 
this condition. The latter sentence means that if $Z \longrightarrow X$ 
factors as $Z \longrightarrow Z' \longrightarrow X$, where $Z'$ is Galois and 
connected, then actually $Z = Z'$.
\end{definition}

We now provide a characterization of Galois closures of \'etale maps, showing 
that the Galois closure always exists.

\begin{definition}
	\label{def:galois_scheme}
	Let $f \colon X \longrightarrow Y$ be a finite \'etale map of degree~$d$ 
	between two irreducible smooth varieties. We define 
	the \emph{Galois scheme} (see \cite[Section~$3$]{Vakil2006}) of~$f$ as
	the scheme $\GS(f)$ such that for any extension~$K$ of the ground field 
	of~$X$ and~$Y$, we have
	\[
	\GS(f)(K) =
	\bigl\{
	(x_1, \dotsc, x_d) \in X^d(K) 
	\, : \, 
	f(x_1) = \dotsb = f(x_d), \; 
	x_i \neq x_j \text{ for all } i \neq j
	\bigr\}.
	\]
\end{definition}

Notice that the Galois scheme is the fiber product of~$d$ copies of the 
map~$f$ minus $\binom{d}{2}$ small diagonals. Because of this, and since $f$ is a finite map, 
we have
\begin{equation}
\label{eq:dimension_galois_scheme}
\dim \GS(f) \, = \, \dim X \, = \, \dim Y.
\end{equation}
There is an induced map $F 
\colon \GS(f) \longrightarrow Y$, sending $(x_1, \dotsc, x_d)$ to~$f(x_1)$, 
which is dominant, and each point $y \in Y$ has $d!$ preimages.
Notice that $\GS(f)$ is Galois over~$Y$.

\begin{proposition}
	Let $f \colon X \longrightarrow Y$ be a finite connected \'etale map of 
degree~$d$ between two irreducible smooth varieties. Pick a point $y_0 \in Y$ 
and let $(x_1, \dotsc, x_d)$ be a fixed permutation of the (geometric) fiber 
of~$f$ over~$y_0$. Let $Z$ be the connected component of~$\GS(f)$ containing 
$(x_1, \dotsc, x_d)$. Then $F|_{Z} \colon Z \longrightarrow Y$ is a Galois 
closure of~$f \colon X \longrightarrow Y$.
\end{proposition}
\begin{proof}
 We have to prove that $Z$ is Galois over $Y$ and that $Z$ is minimal with 
respect to this property. Since $\GS(f)$ is Galois over~$Y$, also the 
restriction of~$F$ to any of its connected components is so, hence $Z$ is 
Galois over $Y$. Suppose now that we have a factorization $Z \longrightarrow Z' 
\longrightarrow X$ with $Z'$ Galois over~$Y$. This induces a factorization of 
the inclusion $Z \hookrightarrow \GS(f)$ as $Z \longrightarrow Z' 
\hookrightarrow \GS(f)$, and this implies that $Z = Z'$.
\end{proof}

\begin{remark}
 Note that the permutation group~$S_d$ of $d$ elements is a group of 
automorphisms of~$\GS(f)$ over~$Y$ acting transitively on the fibers of~$F$. Hence the 
stabilizer of~$Z$ under this group is a group of automorphisms of~$F|_{Z}$ 
acting transitively on the fibers (which shows that $Z$ is Galois over~$Y$). It follows that 
the number of irreducible components of the Galois scheme~$\GS(f)$ coincides 
with the number of cosets of this stabilizer.
\end{remark}

We notice that, if we consider the \'etale map of varieties induced by a separable 
extension of fields $K \subseteq L$, then the spectrum of a Galois closure (in the 
field sense) of $K \subseteq L$ satisfies the universal property of the Galois 
closure (in the map sense, namely as in Definition~\ref{def:galois_closure}).

\begin{definition}
Let $f \colon X \longrightarrow Y$ be a finite \'etale morphism 
between irreducible smooth varieties. Since $f$ is dominant, it determines a 
field extension $K(Y) \hookrightarrow K(X)$. We define the \emph{Galois 
group} $\Gal(f)$ of~$f$ to be the Galois group of the extension $K(Y) 
\hookrightarrow E$, where $E$ is a Galois closure (see 
\cite[Remark~4.77]{Rowen2006}) of $K(Y) \hookrightarrow K(X)$.
\end{definition}

\begin{proposition}
\label{prop:galois_isomorphic}
For every finite \'etale morphism $f \colon X \longrightarrow Y$ of smooth 
irreducible varieties the Galois group of~$f$ is the stabilizer of the Galois 
closure $Z$ of~$f$ in the Galois scheme~$\GS(f)$.
\end{proposition}
\begin{proof}
 The proof follows if we can show that the base change of a Galois closure 
is still a Galois closure if it is connected. In fact, if this is true, given a 
Galois closure $Z \longrightarrow Y$ of~$f$, we can consider its base change at 
the generic point of~$Y$. The base change of $Z \longrightarrow X 
\longrightarrow Y$ under this map is $\mathrm{Spec}\, K(Z) \longrightarrow 
\mathrm{Spec}\, K(X) \longrightarrow \mathrm{Spec}\, K(Y)$. We then know that $K(Z)$ 
is a Galois closure of $K(Y) \subseteq K(X)$, and so $\Gal \bigl(K(Z) / K(Y) 
\bigl)$ is~$\Gal(f)$. However, $\Gal \bigl(K(Z) / K(Y) \bigl)$ coincides with 
the stabilizer of~$Z$ in~$\GS(f)$, because base change preserves the group of 
automorphism lying over the base and permuting the fibers (which makes the 
corresponding map Galois over the base).

Hence, we need to show that the base change of a Galois closure 
is still a Galois closure if it is connected. Suppose that $g \colon W 
\longrightarrow Y$ is a morphism. Then the map $f' \colon X \times_Y W 
\longrightarrow W$ is \'etale by \cite[Lemma 38.34.4]{stacks-project}. Assume 
that $X \times_Y W$ is connected; then by base change $Z \times_Y W$ is a union 
of connected components of~$\GS(f')$. Therefore $Z \times_Y W$ is a Galois 
closure of~$f'$ when it is itself connected.
\end{proof}

Now we cast the notions defined so far into the framework of Galois schemes of 
morphisms (Corollary~\ref{cor:galois_curve_as_map}). After that, we recall the 
notion of simple tangency for a curve and highlight its consequences on Galois 
groups. 

\begin{definition}
\label{def:etale_locus}
For an absolutely irreducible curve~$C \subseteq \p^2_{\F_q}$ of degree~$d$,
define $\mcal{V}_C$ to be the set of points in~$X_0 = \check{\p}^2_{\F_q}$ such that 
the restriction of the map $\pi \colon X_1 \longrightarrow X_0$ 
from Definition~\ref{def:projection} to $\mcal{U}_{C} := 
\pi^{-1}\bigl(\mcal{V}_C\bigr)$ is \'etale.
\end{definition}

\begin{remark}
\label{rem:etale}
Notice that the set $\mcal{V}_C$ is open and non-empty. In fact, since the map 
$\pi \colon X_1 \longrightarrow X_0$ is separable, the general point 
of~$X_0$ belongs to~$\mcal{V}_C$. Moreover, by 
Lemma~\ref{lemma:projection_quasi_finite} the map~$\pi$ is quasi-finite, and 
since both~$X_0$ and~$X_1$ are projective varieties, it is finite, hence 
closed. The locus of point in~$X_1$ where $\pi$ is ramified is closed (since it 
is locally defined by the vanishing of the minors of a Jacobian matrix), so its 
image under~$\pi$ is closed, too. Therefore the locus in~$X_0$ over which $\pi$ 
is unramified is open and non-empty. It is then enough to ensure that
$\pi \colon \mcal{U}_C \longrightarrow \mcal{V}_C$ is flat. Now, the locus in 
the domain where a map is flat is open (see 
\cite[\href{http://stacks.math.columbia.edu/tag/0398}{Tag 0398}, Theorem 
36.15.1,]{stacks-project}), and flat maps are open morphisms (see 
\cite[\href{http://stacks.math.columbia.edu/tag/01U2}{Tag 01U2}, Lemma 
28.24.9]{stacks-project}), so this shows that $\mcal{V}_C$ is open. The fact 
that $\mcal{V}_C$ is non-empty is ensured by the generic flatness result (see 
\cite[\href{http://stacks.math.columbia.edu/tag/0529}{Tag 0529}, Proposition 
28.26.1]{stacks-project}).
\end{remark}

\begin{lemma}
\label{lemma:restriction_etale}
Let $C$ be an absolutely irreducible curve of degree~$d$ over~$\F_q$. 
Then the restriction to $\mcal{U}_{C} := \pi^{-1}\bigl(\mcal{V}_C\bigr)$ of the 
map $\pi \colon X_1 \longrightarrow X_0$ from Definition~\ref{def:projection} is 
a finite separable dominant \'etale morphism between smooth absolutely 
irreducible varieties.
\end{lemma}
\begin{proof}
We know from Section~\ref{preliminaries:galois} that both $X_0$ and $X_1$ are 
smooth and absolutely irreducible. Since $\mcal{V}_C$ and $\mcal{U}_C$ are open 
and non-empty, the same is true for them. Moreover, $\pi$ is a quasi-finite 
separable dominant morphism between projective varieties 
(Lemma~\ref{lemma:projection_quasi_finite}) and so it is finite. Hence, the same 
holds for its restriction~$\pi_{|_{\mcal{U}_C}}$. By Remark~\ref{rem:etale}, the 
map is \'etale, and this concludes the proof.
\end{proof}

By unravelling the definitions, in the light of 
Lemma~\ref{lemma:restriction_etale} we obtain:

\begin{corollary}
\label{cor:galois_curve_as_map}
For an absolutely irreducible projective plane curve~$C$ over~$\F_q$, we have
$\Gal(C) \cong \Gal \bigl( \pi_{|_{\mcal{U}_C}} \bigr)$.
\end{corollary}

The interpretation of the Galois group of a curve provided by 
Corollary~\ref{cor:galois_curve_as_map} allows to use 
Proposition~\ref{prop:galois_isomorphic} and hence to deduce the irreducibility 
of the Galois scheme when the Galois group is the full symmetric group.

\begin{corollary}
Suppose that $C$ is an absolutely irreducible curve in~$\p^2_{\F_{q}}$ of 
degree~$d$ with simple tangency. Then, the Galois group~$\Gal \bigl( 
\pi_{|_{\mcal{U}_C}} \bigr)$ is the full symmetric group, and so the Galois 
scheme $\GS \bigl( \pi_{|_{\mcal{U}_C}} \bigr)$ is  irreducible.
\end{corollary}
\begin{proof}
This follows from Corollary~\ref{cor:galois_curve_as_map} and 
Proposition~\ref{prop:galois_isomorphic}.
\end{proof}

\section{Probabilities of incidence}
\label{probabilities}

In this section we define probabilities of intersection between a random line
and a given curve in the projective plane over a finite field 
(Definition~\ref{def:density}). We then prove the main result of our 
paper, namely Theorems~\ref{thm:existence} and~\ref{thm:formula}, by showing 
that its counterpart for morphisms hold (Theorems~\ref{thm:existence_maps} 
and~\ref{thm:formula_maps}). We will re-prove these results in 
Section~\ref{chebotarev} by using the Chebotarev density theorem.

\begin{definition}[Probabilities of intersection]
\label{def:density}
Let $q$ be a prime power and let $C$ be a plane projective absolutely 
irreducible curve of degree~$d$ over~$\F_q$. 
For every $N \in \N$ and for every $k \in \{0, \dotsc, d\}$, the \emph{$k$-th 
probability of intersection}~$p_k^N(C)$ of lines with~$C$ over~$\F_{q^N}$ is 
\[
p_k^N(C) := 
\frac{\Bigl| 
\bigl\{
\text{lines } \ell \subseteq \p^2_{\F_{q^N}} \, : \, 
|\ell(\F_{q^N}) \cap C(\F_{q^N}) | = k 
\bigr\} 
\Bigr|}{q^{2N}+q^N+1} \,.
\]
Notice that $q^{2N}+q^N+1$ is the number of lines in~$\p^2_{\F_{q^N}}$.
\end{definition}

The aim of this paper is to prove that the limit as~$N$ goes to infinity of the 
quantities~$p_k^N(C)$ exists for every~$k$, and to give a formula for these 
limits, provided that some conditions on the curve~$C$ are fulfilled.

The following result is a direct consequence of Definitions~\ref{def:density} 
and~\ref{def:projection}.

\begin{lemma}
Let $C$ be a plane projective absolutely irreducible curve of 
degree~$d$ over~$\F_q$. 
For every $k \in \{0, \dotsc, d\}$ we have
\[
p_k^N(C) = 
\frac{\Bigl| 
	\bigl\{ 
	[\ell] \in \check{\p}^2 (\F_{q^N}) \, : \, 
	| \pi^{-1}([\ell])(\F_{q^N})| = k 
	\bigr\} 
	\Bigr|}{q^{2N}+q^N+1} \,.
\]
\end{lemma}

Via Lemma~\ref{lemma:reduction_etale} and 
Definition~\ref{def:probabilities_map} we reduce the 
problem of computing intersection probabilities for curves to the 
analogous problem for morphisms.

\begin{lemma}
\label{lemma:reduction_etale}
Let $C$ be a plane projective absolutely irreducible curve of degree~$d$ 
over~$\F_q$. Let $\mcal{V}_C \subseteq \check{\p}^2_{\F_q}$ be 
as in Definition~\ref{def:etale_locus}. For every $N \in \N$ and for every $k 
\in \{0, \dotsc, d\}$, define
\[
\tilde{p}_k^N(C) :=
\frac{\Bigl| 
	\bigl\{ [\ell] \in \mcal{V}_C(\F_{q^N})  \, : \, 
	| \pi^{-1}([\ell]) \bigl( \F_{q^N} \bigr)| = k \bigr\} 
	\Bigr|}{|\mcal{V}_C(\F_{q^N})|} \,.
\]
Then $\displaystyle \lim_{N \to \infty} p_k^N(C)$ exists if and only if 
$\displaystyle \lim_{N \to \infty} \tilde{p}_k^N(C)$ exists, in which case 
the two numbers coincide.
\end{lemma}
\begin{proof}
	It is enough to show that the probability for a 
	point to lie in $\p^2(\F_{q^N}) \setminus \mcal{V}_C(\F_{q^N})$ goes to zero 
as 
	$N$ goes to infinity. This is a consequence of the Lang-Weil bound 
	(Theorem~\ref{thm:lang_weil}). In fact, since the complement of $\mcal{V}_C$ has dimension at most~$1$:
	\[
	\frac{\bigl| \p^2_{\F_{q^N}}(\F_{q^N}) \setminus \mcal{V}_C(\F_{q^N}) 
\bigr|}{q^{2N}+q^N+1} 
	\, \sim \,
	\frac{\bigl( c + O(q^{-N/2}) \bigr) \, q^N}{q^{2N}} \to 0 \, ,
	\]
	where the constant~$c$ is the number of irreducible components of the complement of~$\mcal{V}_C$.
\end{proof}

\begin{definition}
	\label{def:probabilities_map}
	Let $f \colon X \longrightarrow Y$ be a finite \'etale morphism of degree~$d$, 
where $q$ is a prime power, between smooth irreducible varieties over~$\F_q$. 
For every $N \in \N$ and for every $k \in \{0, \dotsc, d\}$, we define the 
\emph{$k$-th preimage probability} $p_k^N(f)$ to be
	\[
	p_k^N(f) := 
	\frac{\Bigl| 
		\bigl\{ 
		y \in Y(\F_{q^N}) \, : \, |f^{-1}(y)(\F_{q^N})| = k 
		\bigr\} 
		\Bigr|}{|Y(\F_{q^N})|} \,.
	\]
\end{definition}

Notice that if $C$ is an absolutely irreducible algebraic plane curve of 
degree~$d$, then for every $N \in \N$ and for every $k \in \{0, \dotsc, d\}$ we 
have $\tilde{p}_k^N(C) = p_k^N\bigl(\pi_{|_{\mcal{U}_C}}\bigr)$. Hence, by 
Lemma~\ref{lemma:reduction_etale}, in order to show the existence of the limits 
of $k$-th probabilities of intersections for a curve, it is enough to show the 
existence of $k$-th preimage probabilities for morphisms over~$\F_q$.

\begin{theorem}
	\label{thm:existence_maps}
	Let $f \colon X \longrightarrow Y$ be a finite \'etale morphism of 
degree~$d$, where $q$ is a prime power, between smooth irreducible varieties 
over~$\F_q$. Then for every 
	$k \in \{0, \dotsc, d\}$ the limit as $N$ goes to infinity of the sequence 
	$\bigl( p_k^N(f) \bigr)_{N \in \N}$ exists.
\end{theorem}
\begin{proof}
	We generalize the construction of the Galois scheme of the morphism $f$.
	For every $k \in \{0, \dotsc, d\}$, define $G_k(f)$ to be the scheme
	such that for every extension~$K$ of~$\F_q$, we have
	\begin{align*}
	G_k(f)(K) := 
	\bigl\{ 
	(x_1, \dotsc, x_k) \in X^k(K)
	\, &: \, 
	f(x_1) = \dotsc = f(x_k), \\
	&\phantom{: \ \ } x_i \neq x_j \text{ for all } i \neq j
	\bigr\}.
	\end{align*}
	In particular $G_d(f) = \GS(f)$.
	As we showed for the Galois scheme, see 
	Equation~\eqref{eq:dimension_galois_scheme}, for every~$k$ the 
	variety~$G_k(f)$ has the same dimension of~$X$ and~$Y$.
	There is a natural finite morphism $F_k \colon G_k(f) \longrightarrow Y$, the 
	fiber product of~$f$ with itself $k$ times. A general $\overline{\mathbb{F}}_q$-valued point of~$Y$ has $d(d-1) \dotsb (d-k+1)$ preimages under the map~$F_k$. The main 
	idea of the proof is to compute, in two different ways, the expected 
	cardinality~$\mu_k^N(f)$ of the set of $\mathbb{F}_{q^N}$-rational points of the fiber~$F_k^{-1}(y)$, where $y$ 
is a uniformly distributed random element in~$Y(\mathbb{F}_{q^N})$. On one hand,
	\[
	\mu_k^N(f) = \frac{\bigl| G_k(f)(\F_{q^N}) \bigr|}{\bigl| Y(\F_{q^N}) \bigr|} 
\,.
	\]
	On the other hand, we can express~$\mu_k^N(f)$ in terms of the preimage 
	probabilities:
	\begin{equation}
	\label{eq:expectation}
	\mu_k^N(f) = \sum_{s = k}^{d} s (s-1) \dotsb (s-k+1) \, p_s^N(f) \,.
	\end{equation}
	In matrix form:
	\begin{equation}
	\label{eq:matrix_form}
	\begin{pmatrix}
	\mu_0^N(f) \\ \vdots \\ \mu_d^N(f)
	\end{pmatrix}
	=
	\begin{pmatrix}
	1 & \ast & \cdots & & \cdots & \ast\\
	0 & 1 & \ast & & & \vdots \\
	\vdots & & \ddots & & & \vdots\\
	0 & \cdots & 0 & k! & \ast & \ast\\
	\vdots & & & & \ddots & \vdots \\
	0 & \cdots & & & \cdots & d!
	\end{pmatrix}
	\begin{pmatrix}
	p_0^N(f) \\ \vdots \\ p_d^N(f)
	\end{pmatrix}.
	\end{equation}
	Since the matrix in Equation~\eqref{eq:matrix_form} has non-zero determinant, 
	we can write 
	\begin{equation}
	\label{eq:probabilities_expectation}
	p_k^N(f) = \sum_{s = 0}^d \alpha_{k,s} \, \mu_s^N(f)
	\end{equation}
	for some numbers~$(\alpha_{k,s})_{k,s}$. Using the Lang-Weil bound on 
	Equation~\eqref{eq:expectation}, we have
	\begin{equation}
	\label{eq:limit}
	\mu_k^N(f) \, \sim \, \frac{\delta_k \, q^{N \cdot \dim G_k(f)}}{q^{N \cdot \dim 
			Y}} 
	\qquad \text{as } N \to \infty \, ,
	\end{equation}
	where $\delta_k$ is the number of irreducible components of 
	$G_k(f)$ that are absolutely irreducible.
	Since $\dim G_k(f) = \dim Y$, we conclude that the 
	limit in Equation~\eqref{eq:limit} exists, and so by 
	Equation~\eqref{eq:probabilities_expectation} also $\displaystyle \lim_{N \to 
		\infty} p_k^N(f)$ exists.
\end{proof}

\begin{remark} 
	Theorem~\ref{thm:existence} holds. In fact, the map $\pi_{|_{\mcal{U}_C}}$ 
satisfies the hypotheses of Theorem~\ref{thm:existence_maps}, so the numbers 
$p_k \bigl( \pi_{|_{\mcal{U}_C}} \bigr)$ exist, and we have already proved that 
this implies that the limits~$p_k(C)$ exist.
\end{remark}

\begin{theorem}
	\label{thm:formula_maps}
	Let $f \colon X \longrightarrow Y$ be a finite \'etale morphism of 
degree~$d$, where $q$ is a prime power, between smooth irreducible varieties 
over~$\F_q$. Suppose that $\Gal(f)$ is the full symmetric group~$S_d$. Then for 
every $k \in \{0, \dotsc, d\}$ we have 
	\[
	p_k(f) = \sum_{s=k}^d \frac{(-1)^{k+s}}{s!} \binom{s}{k}.
	\]
	In particular, $p_{d-1}(f) = 0$ and $p_d(f) = 1 / d!$.
\end{theorem}
\begin{proof}
	Since $\Gal(f)$ is the full symmetric group, the Galois 
	scheme~$\GS(f)$ is absolutely irreducible. Hence, using the notation of the 
	proof of Theorem~\ref{thm:existence_maps}, for all $k \in \{0, 
	\dotsc, d\}$ we have
	\begin{equation}
	\label{eq:expected_preimage}
	\lim_{N \to \infty} \mu_k^N(f) = 
	\lim_{N \to \infty} \frac{q^{N \cdot \dim G_k(f)}}{q^{N \cdot \dim Y}} = 1.
	\end{equation}
	In fact, every variety~$G_k(f)$ is an image (under a projection) of $\GS(f) = 
	G_d(f)$, thus is absolutely irreducible and so 
	Equation~\eqref{eq:expected_preimage} follows from Equation~\eqref{eq:limit}. 
	Again using the notation as in Theorem~\ref{thm:existence_maps}, we get
	\begin{equation}
	\label{eq:formula_alpha}
	\lim_{N \to \infty} p_k^N(f) = \sum_{s = 0}^{d} \alpha_{k,s}.
	\end{equation}
	Therefore, the statement is proved once we are able to explicitly compute the 
	coefficients~$(\alpha_{k,s})_{k,s}$. Recall that $\alpha_{k,s}$ is the 
	$(k,s)$-entry of the inverse of the matrix~$M_d$  
	appearing in Equation~\eqref{eq:matrix_form}. A direct inspection of the 
	matrices~$M_d$ shows that they admit the following structure:
	\[
	M_d = 
	\left(
	\begin{array}{ccc|c}
	\multicolumn{3}{c|}{\multirow{3}{*}{$M_{d-1}$}} & 1 \\
	&&& \vdots \\
	&&& d! / 1! \\ \cline{1-3}
	0 & \cdots & 0 & d! / 0!
	\end{array}
	\right) \,.
	\]
	A direct computation shows that 
	\[
	M_d^{-1} = 
	\left(
	\begin{array}{ccc|c}
	\multicolumn{3}{c|}{\multirow{3}{*}{$M_{d-1}^{-1}$}} & \frac{(-1)^d}{d!} \cdot
	\binom{d}{0} \\
	&&& \vdots \\
	&&& \frac{(-1)}{d!} \cdot \binom{d}{d-1} \\ \cline{1-3}
	0 & \cdots & 0 & \frac{1}{d!} \cdot \binom{d}{d}
	\end{array}
	\right) \,.
	\]
	Hence 
	\[
	\alpha_{k,s} = \frac{(-1)^{k+s}}{s!} \binom{s}{k} 
	\quad \text{ for all } k,s \in \{0, \dotsc, d\}.
	\]
	It follows from Equation~\eqref{eq:formula_alpha} that for all $k \in \{0, 
	\dotsc, d\}$,
	\[
	p_k(f) = \sum_{s=0}^d \frac{(-1)^{k+s}}{s!} \binom{s}{k} 
	= \sum_{s=k}^d \frac{(-1)^{k+s}}{s!} \binom{s}{k}
	\]
	and so the statement is proved.
\end{proof}

As a consequence of Proposition~\ref{prop:symmetric_galois} and 
Theorem~\ref{thm:formula_maps}, Theorem~\ref{thm:formula} holds.

\section{Probabilities of intersection via the Chebotarev theorem}
\label{chebotarev}

In this section, we show how to use an effective version of the Chebotarev 
density theorem for function fields as exposed in 
\cite[Appendix~A]{Andrade2015}---and used in~\cite{Bary-Soroker2012} 
and~\cite{Diem2012} to prove the results reported in 
the Introduction---to show Theorems~\ref{thm:existence} and~\ref{thm:formula}. 
We recall the setting and the results of the paper~\cite{Andrade2015},
and specialize the Chebotarev theorem to our 
case. We refer to the cited appendix for the proofs of the claims we make in 
this section regarding the objects introduced to state the Chebotarev theorem 
(Theorem~\ref{thm:chebotarev}). 

We start by considering an integrally closed finitely generated $\F_q$-algebra 
$R$ and a monic polynomial~$\mcal{F} \in R[T]$ such that the discriminant 
of~$\mcal{F}$ is invertible in~$R$. In our case, we take~$R$ to be the 
$\F_q$-algebra
\[
 R := \frac{\F_q[a,b,u]}{\Disc_x \bigl(F(x, ax+b) \bigr) \cdot u - 1} \cong 
 \F_q[a,b]_{(f)} \quad \text{with } f := \Disc_x \bigl(F(x, ax+b) \bigr),
\]
where the last ring is the localization of the polynomial ring~$\F_q[a,b]$ at 
the element~$f$. In geometric terms, $R$ is the coordinate ring of the open 
subset of the dual projective plane parametrizing lines in the plane that 
intersect the curve $\{F = 0\}$ in $d$ distinct points over the algebraic 
closure of~$\F_q$. We then take the polynomial~$\mcal{F}$ to be $F(T, aT + b)$. 
Then by construction, its discriminant is invertible in~$R$.

Starting from $R$ and~$\mcal{F}$, we consider~$K$, the quotient field of~$R$, 
and we define $L$ to be the splitting field of~$\mcal{F}$ over~$K$. In other 
words, if $\{y_1, \dotsc, y_d\}$ are the roots of~$\mcal{F}$, we set $L := 
K(y_1, \dotsc, y_d)$. In our situation, we have
\[
 L = \frac{K[t_1, \dotsc, t_d]}{\bigl( F(t_i, a t_i + b) \text{ for } i \in 
\{1, \dotsc, d \} \bigr)} \,.
\]
Then we define $S$ to be the integral closure of~$R$ in~$L$, namely $S = R[y_1, 
\dotsc, y_d]$. Geometrically, $S$ is the coordinate ring of an open subset of 
the unique variety $X_d \subset C^d \times \check{\p}^2_{\F_q}$ such that for 
every extension~$M$ of~$\F_q$ we have
\[
 X_d(M) = \bigl\{ (x_1, \dotsc, x_d, [\ell]) \in C^d(M) \times 
\check{\p}^2(M) \, \colon \, x_i \in \ell \bigr\} \,.
\]

The strategy we adopt to compute probabilities of intersections is the 
following: our goal is to count the number of lines~$\ell$ in~$\p^2$ such that 
the intersection $\ell(\F_{q^N}) \cap C(\F_{q^N})$ is constituted of exactly~$k$ points, 
and we interpret this as the number of lines such that the 
univariate polynomial~$F_{|{\ell}}$ has exactly~$k$ linear factors 
over~$\F_{q^N}$. Notice that every univariate 
polynomial~$H$ of degree~$d$ over~$\F_q^{N}$ determines a partition~$\pi_H$ 
of~$d$, namely a tuple $\pi_H = (\alpha_1, \dotsc, \alpha_s)$ such that 
$\alpha_1 + \dotsb + \alpha_s = d$ and $\alpha_1 \leq \dotsb \leq \alpha_s$. 
Such partition is obtained by factoring~$H$ over~$\F_{q^N}$ into irreducible 
factors $H_1, \dotsc, H_s$ and then setting $\alpha_i = \deg(H_i)$. Then, the 
number of lines we are interested in can be computed as the sum, over the set 
of partitions~$\pi$ of~$d$ with exactly~$k$ ones, of the number of lines~$\ell$ 
such that the partition associated to~$F_{|\ell}$ is~$\pi$. The Chebotarev theorem 
provides a formula for the probability for a line to determine a given 
partition.

We set $G$ to be the Galois group of the field extension $K \subseteq L$. By 
definition, this coincides with the Galois group of the curve~$C$ as in 
Definition~\ref{def:galois_curve}.
Notice that, in our situation, the intersection $L \cap \F$, where $\F$ 
is an algebraic closure of~$\F_q$, coincides with $\F_q$. This implies that the 
subgroup
\[
 G_0 := \bigl\{ g \in G \, \colon \, g_{|_{\F_q}}(x) = x \text{ for all } x \in 
\F_q \bigr\}
\]
coincides with~$G$. Similarly, if for every $\nu \geq 1$ we set
\[
 G_\nu := \bigl\{ g \in G \, \colon \, g_{|_{\F_q}}(x) = x^{q^\nu} \text{ for 
all } x \in \F_q \bigr\} \, ,
\]
then $G_{\nu}$, which in general is a coset of~$G_0$ in~$G$, coincides with~$G$.

As one can see from the definition of~$X_d$, its points are intimately related 
to the probabilities we are interested in. From an algebraic point of view 
(see~\cite[Section~II.6]{Mumford1999}) these points correspond to 
$\F_q$-homomorphisms from~$S$ to~$\F$. Moreover, a homomorphism $\Phi \in 
\Hom_{\F_q}(S, \F)$ such that $\Phi(R) = \F_{q^\nu}$ corresponds to a point 
$(x_1, \dotsc, x_d, [\ell])$ in~$X_d$ such that the line $\ell$ is defined 
over~$\F_{q^\nu}$. Given such a homomorphism~$\Phi$ there always exists an 
element in~$G$, called the \emph{Frobenius element} and denoted 
$\left[ \frac{S / R}{\Phi} \right]$ such that the following diagram is 
commutative:
\begin{equation}
\label{eq:frobenius}
 \xymatrix@C=1.5cm{
 S \ar[r]^{\left[ \frac{S / R}{\Phi} \right]} \ar[d]_{\Phi} & S \ar[d]^{\Phi} \\
 \F \ar[r]^{\alpha \mapsto \alpha^{q^{\nu}}} & \F 
 }
\end{equation}
In other words, we have the relation
\[
 \Phi \left( \left[ \frac{S / R}{\Phi} \right] \, x \right) = 
 \Phi(x)^{q^{\nu}} \, .
\]
One then can show that $\left[ \frac{S / R}{\Phi} \right] \in G_{\nu}$.

If we fix a line in~$\p^2$, namely, if we fix an $\F_q$-homomorphism $\varphi 
\in \Hom_{\F_q}(R, \F)$, we can consider all points in~$X_d$ ``lying over''
this line. In other terms, we can consider all homomorphisms $\Phi \in 
\Hom_{\F_q}(S, \F)$ prolonging~$\varphi$. Their corresponding Frobenius 
elements form one key object in the statement of the Chebotarev theorem. For 
$\varphi \in \Hom_{\F_q}(R, \F)$, we set
\[
 \left( \frac{S / R}{\varphi} \right) := \left\{ \left[ \frac{S / 
R}{\Phi} \right] \, \colon \, \Phi \text{ prolongs } \varphi \right\} \, .
\]
In our setting, since $G_0 = G$ one can show that $\left( \frac{S / R}{\varphi} 
\right)$ is a conjugacy class in~$G$. Now we are ready to state the Chebotarev 
theorem (see \cite[Theorem~A.4]{Andrade2015}): 

\begin{theorem}
\label{thm:chebotarev}
Let $Z \subseteq G$ be a conjugacy class and let $\nu \geq 1$; define
\[
 P_{\nu, Z} := \frac{\left| \bigl\{ \varphi \in \Hom_{\F_q}(R, \F) \text{ 
such that } \varphi(R) = \F_{q^{\nu}} \text{ and } \left( \frac{S / R}{\varphi} 
\right) = Z \bigr\} \right|}{\left| \bigl\{ \varphi \in \Hom_{\F_q}(R, \F) 
\text{ such that } \varphi(R) = \F_{q^{\nu}} \bigr\} \right|} \, .
\]
Then there exists a constant~$\delta$ independent of~$q$ such that, as $q \to 
\infty$,
\[
 P_{\nu, Z} \sim \frac{|Z|}{|G|} + \frac{\delta}{\sqrt{q}} \, .
\]
\end{theorem}
In order to use the Chebotarev theorem for our purposes, we have to understand what 
does the condition $ \left( \frac{S / R}{\varphi} \right) = Z$ correspond to in 
our setting. Suppose that $C$ has simple tangency. Then we know by Proposition 
\ref{prop:symmetric_galois} that $G$ is the symmetric group~$S_d$. Notice that 
to every conjugacy class~$Z$ of~$S_d$ we can associate a partition~$\pi_Z$ of 
$d$, obtained from the cycle structure of permutations belonging to~$Z$. On the 
other hand, given a line $\ell = \{ y = ax + b\}$, we can consider the 
restriction of the equation~$F$ of~$C$ to~$\ell$, namely the univariate 
polynomial $F_\ell = F(x, ax + b)$. This polynomial defines a 
partition~$\pi_\ell$ of $d$ by considering its factorization over~$\F_{q^\nu}$: 
the partition~$\pi_\ell$ has as many $1$ as the linear factors of~$F_\ell$, as 
many~$2$ as the quadratic factors of~$F_\ell$, and so on. 

\begin{lemma} 
If $Z \subseteq S_d$ is a conjugacy class of permutations, then the set  
\[
I_{\nu, Z} := \bigl\{ \varphi \in \Hom_{\F_q}(R, \F) \text{ 
	such that } \varphi(R) = \F_{q^{\nu}} \text{ and } \left( \frac{S / 
R}{\varphi} 
\right) = Z \bigr\}
\]
corresponds to the set of lines in~$\p^2_{\F_{q^{\nu}}}$ such that 
$\pi_\ell = \pi_Z$.
\end{lemma}
\begin{proof}
Let $\varphi \in I_{\nu, Z}$ and let $\Phi \in \Hom(S, \F)$ be a homomorphism 
prolonging~$\varphi$. Let $\ell = \{ y = \bar{a} x + \bar{b} \}$ be the line in 
$\p^2_{\F_{q^{\nu}}}$ corresponding to $\varphi$. Then from the explicit 
description of~$K$ and~$L$ we provided at the beginning of the section, it 
follows that $M := \Phi(S)$ is a splitting field of the polynomial $F_\ell = 
F(x, \bar{a} x + \bar{b})$. By definition of the Frobenius element, we have the 
commutative diagram
\[
 \xymatrix@C=1.5cm{ 
   L = \frac{K[t_1, \dotsc, t_d]}{\bigl( 
     F(t_i, a t_i + b) \text{ for } i \in \{1, \dotsc, d \} 
    \bigr)} \ar[d] \ar[r]^{\left[ \frac{S / R}{\Phi} \right]} & 
   L = \frac{K[t_1, \dotsc, t_d]}{\bigl( 
     F(t_i, a t_i + b) \text{ for } i \in \{1, \dotsc, d \} 
    \bigr)} \ar[d] \\
   M = \frac{K[u_1, \dotsc, u_d]}{\bigl( 
     F(u_i, \bar{a} u_i + \bar{b}) \text{ for } i \in \{1, \dotsc, d \} 
    \bigr)} \ar[r]^{\alpha \mapsto \alpha^{q^\nu}} & 
   M = \frac{K[u_1, \dotsc, u_d]}{\bigl( 
     F(u_i, \bar{a} u_i + \bar{b}) \text{ for } i \in \{1, \dotsc, d \} 
    \bigr)} 
 }
\]
which is just the extension to $L$ of the diagram in 
Equation~\eqref{eq:frobenius}. From the commutativity of this diagram, we see 
that the permutation action of $\left[ \frac{S / R}{\Phi} \right]$ on the 
classes $[t_1], \dotsc, [t_d]$ is the same as the action of the map $\alpha 
\mapsto \alpha^{q^\nu}$ on the classes $[u_1], \dotsc, [u_d]$. Since the $\{ 
[u_i] \}$ are the roots of $F(x, \bar{a} x  + \bar{b})$, and the latter is a 
polynomial with coefficients in $\F_{q^\nu}$, which are hence preserved by the 
map $\alpha \mapsto \alpha^{q^\nu}$, it follows that the structure of factors 
of $F_{\ell}$ over $\F_{q^{nu}}$ is the same as the cycle structure of $\left[ 
\frac{S / R}{\Phi} \right]$. This concludes the proof.  
\end{proof}

As a corollary, we obtain that the set of lines in~$\p^2_{\F_{q^\nu}}$ 
intersecting~$C$ in exactly~$k$ points corresponds to the set
\[
 \bigcup_{Z \text{ has exactly } k \text{ fixed points}} I_{\nu, Z} \,.
\]
The number of permutations having exactly~$k$ fixed points is given by the 
so-called \emph{rencontres numbers}, see~\cite{Riordan2002}. We have hence:
\[
 \left| \bigcup_{Z \text{ has exactly } k \text{ fixed points}} Z 
\right| = d! \sum_{s=k}^d \frac{(-1)^{k+s}}{s!} \binom{s}{k} \,.
\]
Using the Chebotarev theorem we then conclude the proof of 
Theorem~\ref{thm:formula}.

As the reader can see, there is nothing particularly special in considering the 
setting of plane curves. In fact, the concept of simple tangency 
(see Definition~\ref{def:simple_tangency}) is applicable to curves in arbitrary projective 
space: an absolutely irreducible curve~$C$ in~$\p^n$ has simple tangency 
if there exists a hyperplane $H \subseteq \p^n_{\overline{\F}_q}$ 
intersecting~$C$ in $d-1$ smooth points of $C$ such 
that~$H$ intersects $C$ transversely at~$d-2$ points and has intersection 
multiplicity~$2$ at the remaining point. Also the concepts of Galois group of a 
curve and probabilities of intersections generalize similarly by considering 
hyperplanes instead of lines.

The generalized statement for the situation of irreducible curves is the 
following.
\begin{proposition}
\label{prop:space_formulas}
Let $C$ be an absolutely irreducible algebraic curve of degree~$d$ in~$\p^n_{\F_q}$, 
where $q$ is a prime power. Suppose that $C$ has simple 
tangency. Then for every $k \in \{0, \dotsc, d\}$ we have 
\[
p_k(C) = \sum_{s=k}^d \frac{(-1)^{k+s}}{s!} \binom{s}{k}.
\]
In particular, $p_{d-1}(C) = 0$ and $p_d(C) = 1 / d!$.
\end{proposition}
Using Proposition~\ref{prop:space_formulas}, we can compute the probabilities of 
intersection of a given plane curve $C$ with a random plane curve~$E$ of 
degree~$e$. In fact, via the \emph{Veronese map} we can reduce this situation to 
the one of Proposition~\ref{prop:space_formulas}. Let us start by defining the 
probabilities of intersection of a given curve~$C$ with a random curve~$E$ in 
the plane in exactly~$k$ points, for $k \in \{ 0, \dotsc, de\}$:
\begin{gather*}
p_k^N(C, e) := 
\frac{\Bigl| 
	\bigl\{
	\text{curves } E \subseteq \p^2_{\F_{q^N}} \text{ of degree } e \, : \, 
	|E(\F_{q^N}) \cap C(\F_{q^N}) | = k 
	\bigr\} 
	\Bigr|}{q^{\binom{e+2}{2}N} + \dotsb + q^{2N}+q^N+1} \,, \\
 p_k(C,e) := \lim_{N \to \infty} p_k^N(C, e) \quad \text{when the limit exists} \,.
\end{gather*}
Recall now that for every $r \in \N$, the Veronese map of degree~$e$ is an 
algebraic morphism embedding~$\p^r$ into a larger projective space, so that 
hypersurfaces of degree~$e$ get mapped to hyperplane sections of the image of 
the map. In this sense, the Veronese map operates a sort of  ``linearization'' 
of the problem. In the case of~$\p^2$, which is the one that interests us, it 
is given by
\[
 \begin{array}{lccc}
  v_e \colon & 
  \p^2  & \longrightarrow & \p^{\binom{e+2}{2}-1} \\
  & (x:y:z) & \mapsto & \bigl( \{ x^a y^b z^c \}_{a+b+c = e} \bigr)
 \end{array} \;.
\]
The following lemma ensures that if we start from a plane curve that has simple 
tangency and we apply the Veronese map, we obtain a curve that has simple 
tangency.

\begin{lemma}
Let $C$ be a plane curve of degree~$d$ with simple tangency and let $e\in \N$. 
Then the image $\tilde{C}=v_{e}(C)$ of $C$ under the Veronese map of degree~$e$ 
has also simple tangency.
\end{lemma}

\begin{proof}
	Let $\ell_1$ be a line witnessing simple tangency for~$C$. Select lines 
$\ell_2, \dots ,\ell_e$ in~$\p^2$ such that each of them intersects~$C$ in~$d$ 
distinct points and $\ell_i \cap \ell_j \cap C$ is empty for all $i \not= j$. 
Define $E$ as  the zero set of the product $\ell_1 \cdots \ell_e$. The Veronese 
map sends~$E$ to a hyperplane section of the Veronese surface; let $\tilde{H}$ 
be the corresponding hyperplane. Then, by construction, $\tilde{H}$ witnesses 
simple tangency for $\tilde{C}$.
\end{proof}

Since the Veronese map of degree~$e$ defines a bijection between plane curves of 
degree~$e$ and hyperplanes in~$\p^{\binom{e+2}{2}-1}$, determining the 
probabilities~$p_k(C,e)$ of intersection of a given plane absolutely 
irreducible curve~$C$ with a random curve of degree~$e$ in~$k$ points is 
equivalent to compute the corresponding probabilities of intersection of the 
image~$\tilde{C}$ of~$C$ under the Veronese map with hyperplanes. We sum up 
what we obtained in the following:
\begin{proposition}
\label{prop:higher_formulas}
Let $C$ be an absolutely irreducible algebraic curve of degree~$d$ in~$\p^2_{\F_q}$, 
where $q$ is a prime power. Suppose that $C$ has simple tangency. 
Let $e \in \N$ be a natural number. Then for every $k \in \{0, \dotsc, de\}$ we 
have 
\[
 p_k(C,e) = \sum_{s=k}^{de} \frac{(-1)^{k+s}}{s!} \binom{s}{k}.
\]
\end{proposition}

\section{Acknowledgements}
We thank all the anonymous referees for very useful comments and suggestions. We 
also thank Claus Diem who shared with us a better strategy to prove 
the results in Section~\ref{galois}, and gave us several other important suggestions to improve significantly the paper.
We also thank Kaloyan Slavov for several useful suggestions and corrections to our paper.
Matteo Gallet has been supported by the Austrian Science Fund 
(FWF): W1214-N15 Project DK9, P26607, P25652, and P31061, and is supported by 
the Austrian Science Fund (FWF): Erwin Schr\"odinger Fellowship J4253. Mehdi 
Makhul has been supported by the Austrian Science Fund (FWF): W1214-N15 Project 
DK9, and is supported by the Austrian Science Fund (FWF): P30405-N32. Josef 
Schicho is supported by the Austrian Science Fund (FWF): W1214-N15 
Project DK9.

\end{document}